\documentclass[a4paper,leqno,10pt,final]{amsproc2000}
\usepackage{amssymb,latexsym,amsxtra2000,verbatim}
\DeclareMathOperator{\Proj}{Proj}

\DeclareMathOperator{\Ext}{Ext}
\DeclareMathOperator{\Hom}{Hom}
\DeclareMathOperator{\RHom}{RHom}

\DeclareMathOperator{\GKdim}{GKdim}

\DeclareMathOperator{\character}{char}

\DeclareMathOperator{\coh}{coh}
\DeclareMathOperator{\Qch}{Qch}
\DeclareMathOperator{\gr}{gr}
\DeclareMathOperator{\Gr}{Gr}
\DeclareMathOperator{\tors}{tors}
\DeclareMathOperator{\Tors}{Tors}
\DeclareMathOperator{\qgr}{qgr}
\DeclareMathOperator{\QGr}{QGr}

\newcommand{\op}{{\mathrm op}}

\newcommand{\ZZ}{{\mathbb Z}}

\newcommand{\RR}{{\mathbb R}}

\newcommand{\NN}{{\mathbb N}}
\newcommand{\PP}{{\mathbb P}}
\newcommand{\LL}{{\mathcal L}}
\newcommand{\calF}{{\mathcal F}}

\newcommand{\calE}{{\mathcal E}}
\newcommand{\calA}{{\mathcal A}}
\newcommand{\calB}{{\mathcal B}}

\newcommand{\calC}{{\mathcal C}}

\newcommand{\calN}{{\mathcal N}}
\newcommand{\calM}{{\mathcal M}}

\newcommand{\OO}{{\mathcal O}}

\newcommand{\gs}{{\sigma}}

\numberwithin{equation}{section}

\theoremstyle{plain}
\newtheorem{theorem}[equation]{Theorem}

\newtheorem{conjecture}[equation]{Conjecture}

\theoremstyle{definition}

\newtheorem{definition}[equation]{Definition}

\newtheorem{question}[equation]{Question}

\theoremstyle{remark}

\begin{document}

\title{The rings of noncommutative projective geometry}

\author{Dennis S. Keeler}  
     \thanks{
     %Dept. of Mathematics, MIT, Cambridge, MA 02139, \texttt{dskeeler@mit.edu},
     Partially supported by an NSF Postdoctoral Research Fellowship.}
   
\date{} 
\address{ Department of Mathematics, MIT, Cambridge, MA 02139-4307}
% \\ Email: dskeeler@mit.edu   }
\email{dskeeler@mit.edu}
\urladdr{http://www.mit.edu/\~{}dskeeler}

%\keywords{ noetherian  graded rings, noncommutative projective geometry}

\begin{abstract}
In the past 15 years a study of ``noncommutative projective
geometry'' has flourished. By using and generalizing techniques of
commutative projective geometry, one can study certain noncommutative
graded rings and obtain results for which no purely algebraic proof is
known. For instance, noncommutative graded domains of quadratic growth,
or ``noncommutative curves,'' have now been classified by geometric data
and these rings must be noetherian. Rings of cubic growth, or
``noncommutative surfaces,'' are not yet classified, but a rich theory
is currently forming.
In this survey, we describe some of these results and examine the question
of which rings should be included in noncommutative projective
geometry.
\end{abstract}

\maketitle

\section{Introduction}

Noncommutative projective geometry is a new field of mathematics
which studies noncommutative graded rings and their categories of modules
via the techniques of projective algebraic geometry. 
This new ``geometry'' has been able to prove ring-theoretic theorems
(for which no purely algebraic proof is known) and has led to new questions
which a classical ring theorist would not have asked.

This paper will survey some of these new answers and questions. But first,
we would like to say that an excellent survey article \cite{SV}
already exists. In fact, the lecture on which this paper is based was in turn
based on that survey. Thus,  similarities between the two papers must
exist and we express our great debt to the previous
article's authors. We hope that this simplified survey will encourage
the reader to examine \cite{SV}. 
However, we have also added some new material, particularly in Section~\ref{sec:strong-noeth}.
%We will mention some topics not included in \cite{SV}, such as intersection theory and
%abstract Hilbert schemes.
%In $\S$\ref{sec:last}, we will mention some topics not covered the previous survey.

To give this paper an overlying structure, we will give our opinionated
answer to the question
%%% every question, theorem, etc, must be between \begin{} and \end{}
%%% the \label allows the question to be reference later with \ref
\begin{question}\label{qu:main} 
Which  rings should be included in noncommutative projective geometry?
\end{question}
As mentioned before, we want our rings to be graded. More specifically,
we usually assume
\begin{enumerate}
\setcounter{enumi}{-1} 
%%% these manual counter resets are not usually necessary; I just needed extra control
\renewcommand{\theenumi}{R\arabic{enumi}}
\item\label{R0} $k$ is an algebraically closed field, 
$R = \bigoplus_{i=0}^\infty R_i$, $R_0 = k$, $\dim_k R_i < \infty$.
\end{enumerate}
While some work \cite{AZ, Ke-ample-filters} has dealt with
the case where $k$ is a commutative noetherian ring, most work
has assumed that $k$ is a field, often algebraically closed.
A graded ring $R$ is \emph{finitely graded}  when $\dim_k R_i < \infty$
and $R$ is \emph{connected} when $R_0 = k$.
We will take \emph{connected graded} to mean connected and finitely graded.
We will sometimes weaken \eqref{R0} to just finitely graded, but
in this paper,
$k$ is an algebraically closed field unless stated otherwise.

It is clear that \eqref{R0} is much too broad for a successful use
of geometric methods. The most obvious example is that
of polynomials in two noncommuting variables, $k\langle x, y \rangle$.
While this ring is finitely generated over $k$, it is not noetherian
and even contains the subring
$k\langle xy^i \colon i \in \NN \rangle$, which is isomorphic
to polynomials in countably many noncommuting variables.
Since projective geometry relies heavily on the fact that
finitely generated commutative $k$-algebras
are noetherian, we stipulate that
\begin{enumerate}
\setcounter{enumi}{0}
\renewcommand{\theenumi}{R\arabic{enumi}}
\item\label{R1} $R$ is right noetherian.
\end{enumerate}

We shall continue to mark our cumulative answer to Question~\ref{qu:main}
by \eqref{R0}, \eqref{R1}, \dots. If an assumption (R$n$) must be strengthened,
we denote the new assumption (R$n'$). We will base these assumptions
on two guiding principles. First, that these properties hold for 
\emph{twisted homogeneous coordinate rings}, the class of rings
closest to commutative graded rings (in the sense of
noncommutative projective geometry). Second, that these
assumptions lead to useful generalizations of commutative
geometry theorems.

An important specialization of Question~\ref{qu:main} is
\begin{question}\label{qu:qp}
Which rings should generalize $k[x_0, \dots, x_n]$? In other words,
which rings should be called a \emph{noncommutative projective space}
or \emph{quantum $\PP^n$}?
\end{question}
We will denote these extra assumptions as (QP$n$).
We warn the reader that these answers (R$n$) and (QP$n$) are
simply the opinions of the author and should not be
taken as an established answers to Questions~\ref{qu:main} and \ref{qu:qp}.
In particular, the definition of a quantum $\PP^n$ changes
from paper to paper, depending on the needs of each author.

\section{Artin-Zhang Theorem}\label{sec:AZT}

Noncommutative rings have few prime ideals
compared to commutative rings. Hence, defining a projective scheme
as a ringed topological space on the homogeneous primes of $R$
would not be useful. Instead, we look to the Serre Correspondence Theorem
(\textbf{SCT})
for a definition of noncommutative projective scheme.

First, we need some notation for module categories associated
to $R$ or a commutative scheme $X$.
\begin{itemize}
\item $\Gr R$ -- the category of graded right $R$-modules,
\item $\Tors R$ -- the full subcategory generated by all $M \in \Gr R$ with $M_i = 0$ 
for all $i \gg 0$,
\item $\QGr R$ -- the quotient category $\Gr R/\Tors R$,
\item $\gr R, \tors R, \qgr R$ -- the respective
 subcategory of noetherian objects,
\item $\Qch X, \coh X$ -- the category of quasicoherent (resp. coherent) sheaves
on the scheme $X$.
\end{itemize}

The SCT states that
if $R$ is commutative and finitely generated by $R_1$, 
then there is a category equivalence
$\QGr R \cong \Qch(\Proj R)$. % and $\qgr R \cong \coh(\Proj R)$,
%where $\Qch X$ %(resp. $\coh X$)
% is the category of quasicoherent
%(resp. coherent)
% sheaves on a projective scheme $X$.
 Further, under these equivalences,
$\pi R$ corresponds to $\OO_X$, where $\pi R$ is the image of $R$
in $\QGr R$. % or $\qgr R$.
Since the categories $\QGr R$ and $\Qch X$ are determined by
their noetherian objects, the SCT is sometimes stated as
the equivalence $\qgr R \cong \coh X$ \cite[Exercise~II.5.9]{Ha1}.

We wish to generalize the SCT to other categories
so that we may study $\QGr R$ when $R$ is noncommutative.
The SCT is proven by representing $R$ as a homogeneous coordinate
ring in high degree; we will define a ``categorical'' homogeneous coordinate ring.
So let $\calC$ be a locally noetherian
 $k$-linear abelian category (meaning that $\calC$ is generated by its noetherian objects
 and that the objects and
$\Hom$-sets of $\calC$ are $k$-vector spaces),
 let $\OO$ be a noetherian object,
and let $s$ be an autoequivalence. 
We define the homogeneous coordinate ring as
\[
\Gamma(\calC, \OO, s)_{\geq 0} = \bigoplus_{i=0}^\infty \Hom(\OO, s^i \OO)
\]
with multiplication given by the composition of homomorphisms
\begin{multline*}
\Hom(\OO, s^i \OO) \times \Hom(\OO, s^j \OO)  \cong 
\Hom(s^j\OO, s^{i+j} \OO) \times \Hom(\OO, s^j\OO) \\ \to \Hom(\OO,s^{i+j}\OO).
\end{multline*}

There are two categories and two autoequivalences which appear in the SCT.
 When $\calC = \QGr R$ and $\OO = \pi R$,
there is a natural autoequivalence $(+1)$ given by $M \mapsto M(1)$ where
$M(i)$ is the $i$th shift of a graded module $M$, that is, $M(i)_j = M_{i+j}$.
(This is still an autoequivalence when $R$ is noncommutative.)
When  $X$ is a projective scheme with $\calC = \Qch X, \OO = \OO_X$,
and $s = - \otimes \LL$ for an invertible sheaf $\LL$, we have the homogeneous
coordinate ring as it is normally defined.

Of course the proof of the SCT requires  $\LL$ to be ample,
so we need a categorical definition of ampleness.
\begin{definition}\label{def:ample}
Let $\calC$ be a locally noetherian
 $k$-linear abelian category, let $\OO \in \calC$ be a noetherian object,
and let $s\colon \calC \to \calC$ be an autoequivalence.
Write $s^i\calM = \calM(i)$. The pair $(\OO,s)$ is \emph{ample} if
\begin{enumerate}
\renewcommand{\theenumi}{A\arabic{enumi}}
\item\label{A1} For all noetherian $\calM \in \calC$, there exist $\ell_i > 0, i = 1, \dots, p$,
such that there exists an epimorphism 
\[
\bigoplus_{i=1}^p \OO(-\ell_i) \twoheadrightarrow \calM,
\]
\item\label{A2} For all epimorphisms $\calM \twoheadrightarrow \calN$, there exists
$n_0$ such that \[ \Hom(\OO, \calM(n)) \to \Hom(\OO, \calN(n)) \] is an epimorphism
for $n \geq n_0$.
\end{enumerate}
\end{definition}

Finally, we need a technical definition; we allow
$k$ to be any commutative noetherian ring, as we will
need a more general definition for use in Section~\ref{sec:strong-noeth}.
\begin{definition}\label{def:chi}
{\upshape \cite[Definition~3.2, Proposition~3.1(3), 3.11]{AZ}}
Let $k$ be a commutative noetherian ring,
let $R$ be a right noetherian graded $k$-algebra where each $R_i$
is a finite $k$-module. %and let $R_{\geq n} = \bigoplus_{i=n}^\infty R_i$.
We say that $R$ satisfies $\chi_j$
if for all $j' \leq j$, %all $n,d \in \NN$, 
and  all noetherian
right graded $R$-modules $M$,
\[
\bigoplus_{i=-\infty}^\infty \Ext^{j'}_{\Gr R} (R_0, M(i))
\]
is a finite $k$-module.
\begin{comment}
Let $R$ satisfy \eqref{R0} and \eqref{R1}, and
let $R_+ = \oplus_{i>0} R_i$. 
We say that $R$ satisfies $\chi_j$
if for all $i \leq j$ and for all noetherian
right graded $R$-modules $M$, $\dim_k \Ext^i (R/R_+, M) < \infty$
(where $M$ is viewed as an ungraded module when computing $\Ext$).
\end{comment}
If $R$ satisfies $\chi_j$ for all $j \geq 0$, we say that $R$ satisfies $\chi$.
\end{definition}

We now return to assuming that $k$ is an algebraically closed field.
We may now state the Artin-Zhang Theorem, a noncommutative analogue of the SCT.

\begin{theorem}[Artin-Zhang] {\upshape \cite[Theorem~4.5]{AZ}} 
Let $\calC$ be a locally noetherian
 $k$-linear abelian category, let $\OO \in \calC$ be a noetherian object,
and let $s\colon \calC \to \calC$ be an autoequivalence.
Assume $\dim_k \Hom(\calM, \calN) < \infty$
for all noetherian $\calM, \calN \in \calC$, and let $(\OO, s)$ be ample.
Set $B= \Gamma(\calC, \OO, s)_{\geq 0}$. Then
$\calC \cong \QGr B$, and $B$ is right noetherian and has $\chi_1$.

Conversely, if $R$ is finitely graded, right noetherian,
and satisfies $\chi_1$, then we have that
$R \cong \Gamma(\QGr R, \pi R, (+1))_{\geq 0}$ (in sufficiently high degree),
$(\pi R, (+1))$ is ample, and $\dim_k \Hom(\calM, \calN) < \infty$ for all
noetherian $\calM, \calN \in \QGr R$.
\end{theorem}

In light of the above theorem, we add to our list
\begin{enumerate}
\setcounter{enumi}{1}
\renewcommand{\theenumi}{R\arabic{enumi}}
\item\label{R2} $R$ satisfies $\chi_1$.
\end{enumerate}
All commutative noetherian $k$-algebras with \eqref{R0} satisfy $\chi_1$.
Unfortunately, this is not the case for all
noncommutative right noetherian rings $R$. Such
rings can exhibit bad behavior; \cite{StafZ} present a family of rings
which are (right or left) noetherian if and only if $\character k = 0$.

\begin{definition}
Let $R$ satisfy \eqref{R0}--\eqref{R2}. Then
the pair $(\QGr R, \pi R)$ is a \emph{noncommutative projective scheme}.
\end{definition}

One may notice that in Definition~\ref{def:ample}, the property \eqref{A2}
is weaker than the cohomological definition of an ample invertible
sheaf given by the Serre Vanishing Theorem.
For noncommutative projective schemes, this vanishing requires not just
the condition $\chi_1$, but $\chi$.

\begin{theorem}{\upshape \cite[Corollary~7.5]{AZ}}
Let $R$ satisfy \eqref{R0}--\eqref{R2}. %and suppose that
%$\dim_k \Ext^i (\calM, \calN) < \infty$ for all noetherian $\calM, \calN \in \QGr R$. 
Then $R$ satisfies $\chi$ if and only if
for all noetherian $\calM$, 
\begin{enumerate}
\item $\dim_k \Ext^i (\pi R, \calM) < \infty$ for all $i \geq 0$, and
\item there exists $n_0$ such that $\Ext^i(\pi R,\calM(n))=0$
for all $n \geq n_0, i>0$.
\end{enumerate}
\end{theorem}

There are rings which satisfy $\chi_1$ but not $\chi$. We will say a bit more
about one such example in the next section. We end this section with another
definition concerning $\Ext$.
\begin{definition}
Let $\calC$ be a locally noetherian, $k$-linear category and let $\OO$ be
a noetherian object. The pair $(\calC, \OO)$ has \emph{finite cohomological
dimension $n$} if $n$ is minimal with respect to
the property that $\Ext^i(\OO, \calM) = 0$ for all $\calM \in \calC$ and $i > n$.
If $R$ is a ring with \eqref{R0}--\eqref{R2} and
 $(\QGr R, \pi R)$ has finite cohomological dimension $n$, we say that $R$
has finite cohomological dimension $n$.
\end{definition}

\section{Twisted homogeneous coordinate rings}

In this section we will discuss those rings $R$ which are closest to commutative rings
in the sense of the Artin-Zhang Theorem. 
\emph{Twisted homogeneous coordinate rings} are rings of the form
\[
\Gamma(\Qch X, \OO_X, s)_{\geq 0}
\]
for some proper commutative scheme $X$, with $(\OO_X, s)$ ample.
By the Artin-Zhang Theorem, a ring $R$ with
\eqref{R0}--\eqref{R2} and  $(\QGr R, \pi R) \cong (\Qch X, \OO_X)$
for some commutative proper scheme $X$ 
is isomorphic to a twisted homogeneous coordinate ring in sufficiently high degree.

Twisted homogeneous coordinate rings are so called because any autoequivalence
of $\Qch X$ is of the form $\gs_*(- \otimes \LL)$ for some automorphism $\gs$
and invertible sheaf $\LL$ on $X$ \cite[Corollary~6.9]{AZ}, \cite[Proposition~2.15]{AV}.
Knowledge of
automorphisms and invertible sheaves then allows one to study these twisted rings.

But before presenting some results, 
let us define a noncommutative analogue of Krull dimension, the
\emph{Gel'fand-Kirillov dimension} 
\[
\GKdim R = 
\inf \{ \alpha \in \RR \colon \dim_k \sum_{i=0}^n R_i \leq n^\alpha \textrm{ for all }n \gg 0 \}.
\]
If $R$ is commutative, then $\GKdim R$ equals the Krull dimension of $R$ \cite[Theorem~4.5]{KL}.
It turns out that the GK-dimension can be $0$, $1$, any real number $\geq 2$, or infinity
\cite[ Proposition 1.4, Theorems 2.5, 2.9]{KL}.
It is an open question whether 
the GK-dimension of a graded domain must be in $\NN \cup \{\infty\}$.
Further, if $R$ is a noetherian graded domain, must we have $\GKdim R \in \NN$?

If $\GKdim R = 0$, then $R$ is a finite dimensional ring.
If $\GKdim R = 1$, then $R$ is well-understood \cite{SSW-GK1}
and if $R$ is a finitely generated domain, then $R$ is commutative  \cite{SW-GK1}.
One of the most interesting results of noncommutative projective geometry
is the classification of domains with GK-dimension $2$.

\begin{theorem}{\upshape \cite[Theorem~0.2]{AS}} Let $R$ be a domain with \eqref{R0}, 
generated in degree one, and let
 $\GKdim R = 2$. Then there exists a commutative projective curve
$Y$ and autoequivalence $s$ of $\Qch Y$
 such that $R \cong \Gamma(\Qch Y, \OO_Y, s)_{\geq 0}$ (up to finite dimension)
with $(\OO_Y, s)$ ample. Thus $R$ is right noetherian and has $\chi_1$.
(Actually, we will soon see that $R$ is noetherian and has $\chi$.)
\end{theorem}

If $R$ is a domain with $\GKdim R = 2$, but no Veronese subring of $R$ is
 generated in degree one,
then $R$ may be non-noetherian. It is also possible that $R$ does not have $\chi_1$.

Artin and Stafford extended their results to  prime rings $R$ 
and  arbitrary fields $k$ \cite{AS2}.
They obtain $\qgr R \cong \coh \calE$, where $\calE$ is a sheaf of orders on a projective
curve $Y$ inside a central simple $k(Y)$-algebra and $\coh \calE$ is the category
of coherent $\calE$-modules. 
For instance, with $Y = \PP^1$ and $\calE = \OO_Y + M_2(\OO_Y(-1))$, one can generate
an $R$ which is
a noetherian PI algebra which is not finite over its center \cite[0.2]{AS2}.
%The simplest non-trivial example
%would be $\calE = M_2(\OO_X)$ and $R = M_2(A)$ where $A$ is a domain with $\GKdim A = 2$.

We  now present important properties of any twisted homogeneous coordinate ring,
which we believe other rings of noncommutative geometry should have.

\begin{theorem}\label{th:twistedthm} {\upshape \cite{Ke-ample-filters, Ke1}}
 Let $R$ have \eqref{R0}--\eqref{R2}. Suppose that %we have
$(\QGr R,\pi R) \cong (\Qch X, \OO_X)$ for some commutative scheme $X$ which is
proper over $k$. Then $X$ is projective over $k$ and,
\begin{enumerate}
\renewcommand{\theenumi}{R\arabic{enumi}$'$}
\item\label{R1'} $R$ is noetherian,
\item\label{R2'} $R$ and $R^\op$ have $\chi$ (where $R^\op$ is the opposite ring),
\renewcommand{\theenumi}{R\arabic{enumi}}
\item\label{R3} $R$ and $R^\op$ have finite cohomological dimension,
\item\label{R4} some Veronese subring $R^{(d)}$ is generated in degree one, and
\item\label{R5} $\GKdim R$ is an integer.
\end{enumerate}
With the appropriate generalizations of definitions,
if $k$ is only assumed to be a commutative noetherian ring, then \eqref{R1'}--\eqref{R4}
are still true.
\end{theorem}

In the theorem above, it is important that $\pi R$ corresponds to $\OO_X$ under the
category equivalence. There exists a ring $R$ with $(\QGr R, \pi R) \cong (\Qch \PP^1, \calF)$
and $\calF \neq \OO_X$ such that $R$ is right noetherian and has $\chi_1$, but
is not left noetherian and does not have $\chi$ \cite[Proposition~6.13]{AZ}.

We have already mentioned that a ring of GK-dimension $2$ without \eqref{R4} may
behave badly. The strengthened properties \eqref{R1'}--\eqref{R3}
will be sufficient for an analogue of Serre Duality to hold. 
We will also see that \eqref{R1'}--\eqref{R5} are properties
of ``quantum polynomial rings,'' which we now examine. 

\section{Quantum polynomial rings}
\label{sec:qpr}

Since any commutative projective scheme can be embedded in $\PP^n$
for some $n$, many projective geometry propositions are mainly
concerned with projective $n$-space or its homogeneous coordinate
ring $k[x_0, \dots, x_n]$. Thus, we would like to define which rings $R$
should generalize polynomial rings, then call $R$ a \emph{quantum polynomial ring}
and call $(\QGr R, \pi R)$ a
\emph{noncommutative $\PP^n$}. As of yet, there is no consensus on what
this definition should be, so the reader should carefully examine the hypotheses of
each paper dealing with noncommutative $\PP^n$.
We will explain some of the most common hypotheses.

First, since commutative polynomial rings are Gorenstein, we would
like our quantum polynomials to satisfy a noncommutative analogue.

\begin{definition}\label{def:AS-Gorenstein} Let $R$ be connected graded \eqref{R0}.
Then $R$ is  \emph{Artin-Schelter Gorenstein (AS-Gorenstein)} if
\begin{enumerate}
\item $R$ has finite right and left injective dimension $d$, and
\item For some shift $\ell$, $\Ext_{\Gr R}^i(k, R) = \begin{cases} 0 & \textrm{if } i\neq d \\
									k(\ell) & \textrm{if } i = d.
					  \end{cases}$
					  
\end{enumerate}
\end{definition}

Polynomial rings also have finite global dimension, so we make another
definition.

\begin{definition} Let $R$ be connected graded \eqref{R0}.
Then $R$ is an \emph{Artin-Schelter regular
(AS-regular) algebra} of dimension $d$ if
\begin{enumerate}
\item $R$ has right and left global dimension $d$,
\item $\GKdim R < \infty$, and
\item $R$ is AS-Gorenstein with injective dimension $d$.
\end{enumerate}
\end{definition}

\begin{definition}\label{def:qpr}
Let $R$ be connected graded \eqref{R0}. Then $R$ is
a \emph{quantum polynomial ring} of dimension $d$ and $(\QGr R, \pi R)$ is a \emph{noncommutative
$\PP^d$} if
\begin{enumerate}
\renewcommand{\theenumi}{QP\arabic{enumi}}
\item $R$ and $R^\op$ are AS-regular of dimension $d$,
\item If $\ell$ is the shift in Definition~\ref{def:AS-Gorenstein}, then $\ell = d$,
\item\label{QP3-noethdomain} $R$ is a noetherian  domain (so \eqref{R1'} holds),
\item $R$ is generated in degree one (so \eqref{R4} holds), and
\item $R$ has Hilbert series $(1-t)^{-d}$.
\end{enumerate}
\end{definition}

Our definition of quantum polynomial ring is stronger than that of many papers,
but it has some nice consequences.
Noetherian AS-Gorenstein rings satisfy \eqref{R2'} and \eqref{R3} \cite{YZ-Serre}.
Our hypotheses also imply that $\GKdim R$ is not only finite, but
an integer \cite[Theorem~4.1]{J-CM-regular}, so $R$ satisfies
the (strengthened) \eqref{R0}--\eqref{R5}.

Of course this definition is only good if it leads to useful results, and that it does.
For instance, if $R$ is a quantum polynomial ring which is also Cohen-Macaulay (in an appropriate
sense), then the \emph{Castelnuovo-Mumford regularity} of noetherian graded $R$-modules can
be defined and has some of the same properties as commutative regularity.
From this, the degrees of generators of some ideals can be bounded \cite{J-CM-regular}.
There is also an \emph{intersection theory} and \emph{B\'ezout's theorem}
on noncommutative $\PP^n$ (with the condition that $R$ is Auslander-Gorenstein,
which is slightly stronger than AS-Gorenstein) \cite{MS-Bezout}.

We should note that  the condition \eqref{QP3-noethdomain} may be implied by some
of the other (QP$n$). AS-regular algebras of dimension $\leq 3$ have
been classified and they are all noetherian domains \cite{ATV,ATV2,St1,St2}.
(The classification of AS-regular algebras of dimension $3$
was the original question which gave birth to noncommutative projective geometry.)
If the Hilbert series does not equal $(1-t)^{-3}$, then the ring is considered
a \emph{weighted} quantum polynomial ring.

\section{Serre duality}

Perhaps the strongest tool which has carried over from commutative projective
geometry is that of Serre duality. It has been instrumental in developing
the intersection theory of \cite{J-intersection, MS-Bezout}.
It has also been used to
analyze Gorenstein-like conditions; \cite{JZ-gourmet}
gives generalizations of Watanabe's Theorem for AS-regular algebras
and also Stanley's Theorem for ``AS-Cohen-Macaulay'' algebras.

%Classical Serre duality, distinguished object $\OO \in \calC$.
Let us now define classical Serre duality. But first, we need a 
definition of ``properness.''

\begin{definition}\label{def:Ext-finite}
Let $\calC$ be a noetherian $k$-linear abelian category. We call
$\calC$ \emph{$\Ext$-finite} if $\dim_k \Ext^i(\calM, \calN) < \infty$
for all $i$ and all $\calM, \calN \in \calC$. This can
be thought of as saying that $\calC$ is proper over $k$.
\end{definition}

\begin{definition}
Let $\calC$ be $\Ext$-finite with
distinguished object $\OO$. Then 
$(\calC, \OO)$ satisfies \emph{classical Serre duality}
 if there exists $\omega^\circ \in D^b(\calC)$
such that $\RHom(-, \omega^\circ) \cong \RHom(\OO, -)^*$
where ${}^*$ is the $k$-dual.
If $(\calC, \OO) = (\qgr R, \pi R)$ for
a connected graded ring $R$, we say that $R$ satisfies
classical Serre duality.
\end{definition}

Of course this definition will only be useful if classical
Serre duality holds for many $(\qgr R, \pi R)$.
The condition on $\RHom$ forces $R$ to have finite cohomological
dimension, which we denoted as the good property \eqref{R3}.
Fortunately, the addition of some of our other properties 
is sufficient for classical Serre duality.

\begin{theorem}\label{th:SerreDuality} 
{\upshape \cite{YZ-Serre}} Let $R$ be connected graded \eqref{R0}
and noetherian \eqref{R1'}.
Suppose
\begin{enumerate}
\renewcommand{\theenumi}{R\arabic{enumi}$'$}
\setcounter{enumi}{1}
\item $R$ and $R^\op$ satisfy $\chi$, and
\renewcommand{\theenumi}{R\arabic{enumi}}
\item $R$ and $R^\op$ have finite cohomological dimension.
\end{enumerate}
Then $R$ satisfies classical Serre duality.
\end{theorem}

Note in particular that our twisted homogeneous coordinate rings satisfy 
classical Serre duality, as do quantum polynomial rings.
For quantum polynomial rings, the situation is even better, because
when $R$ satisfies (strengthened) \eqref{R0}--\eqref{R3} and is AS-Gorenstein of dimension $d+1$,
then $\omega^\circ$ is the image in $D^b(\qgr R)$ of an $R$-bimodule $\omega(d)$, left and right
free of rank $1$, and the duality statement becomes \cite{YZ-Serre}
\[
\Ext^i_{\qgr R}(-, \omega) \cong \Ext^{d-i}_{\qgr R}(\OO, -)^*,
\]
a familiar form of Serre duality for commutative Gorenstein varieties.

Before ending this section, we mention two other possible generalizations
of Serre duality. First, Theorem~\ref{th:SerreDuality} has a converse
in some sense:
If $R$ is noetherian and connected graded, 
then $(\qgr R,\pi R)$ has a ``balanced dualizing complex'' if and only if $R$ has
\eqref{R2'} and \eqref{R3} \cite{Vdb-Existence,YZ-Serre}.
Second, a more category-theoretic version of Serre duality is given
by a \emph{Serre functor} $F \colon D^b(\calC) \to D^b(\calC)$
which gives natural
isomorphisms $\Hom(\calA,\calB) \cong \Hom(\calB,F\calA)^*$ for 
$\calA, \calB \in D^b(\calC)$ \cite{BK-Serre}.

\section{Strongly noetherian rings and abstract Hilbert schemes}\label{sec:strong-noeth}

Changing the base of a projective scheme from $k$ to another $k$-algebra
is a powerful tool of algebraic geometry. %championed by Grothendieck \cite{EGA}.
There has been some recent work on base change in noncommutative geometry,
which has led in turn to  \emph{abstract Hilbert schemes} which parameterize
certain modules in $\qgr R$.

Of course to have a workable theory of base change, we need our rings to 
behave well under tensor products. We thus define
\begin{definition}
Let $R$ be a $k$-algebra (not necessarily graded). Then $R$
is \emph{strongly right noetherian} if $R \otimes_k A$ is
right noetherian for any commutative noetherian
(not necessarily graded) $k$-algebra $A$.
\end{definition}

\begin{definition}
Let $R$ be a connected graded ring \eqref{R0}. Then $R$
satisfies the \emph{strong $\chi$} condition
if $R \otimes_k A$ satisfies $\chi$ (as an $A$-algebra)
 for any commutative noetherian
(not necessarily graded) $k$-algebra $A$.
\end{definition}

(Our definition of strong $\chi$ is weaker than that in \cite[C6.8]{AZ-Hilbert},
but the statement of Theorem~\ref{th:Hilbert} will still be correct.)

Twisted homogeneous coordinate rings are strongly noetherian 
and they also satisfy strong $\chi$ \cite[Proposition~4.13]{ASZ-Generic}
(while the second fact is not explicitly stated, it follows from
the proof of the first fact). Thus we may wish to amend our list.
\begin{enumerate}
\renewcommand{\theenumi}{R\arabic{enumi}$''$}
\item $R$ is strongly noetherian,
\item $R$ and $R^\op$ satisfy strong $\chi$.
\end{enumerate}

Unfortunately, not all noetherian $k$-algebras are strongly noetherian.
\cite{RS-affine} gave a nongraded example. Recently,
Jordan introduced a family of connected graded rings which are not strongly noetherian
\cite{Jordan-eulerian} and Rogalski showed that these rings are noetherian
\cite{Rog-generic-surface}.

Besides giving a version of generic flatness for noncommutative rings
\cite{ASZ-Generic}, the property of strongly noetherian also allows
for the formation of Hilbert schemes.

\begin{theorem}\label{th:Hilbert} {\upshape \cite[Theorem~E5.1]{AZ-Hilbert}}
Let $R$ be a connected graded, strongly noetherian ring, satisfying strong $\chi$. Let $P$ be
a finitely graded, noetherian $R$-module. The quotients of $\pi P$ in $\qgr R$ with fixed
Hilbert function $h$ are parameterized by a countable union of commutative projective schemes.
\end{theorem}

\begin{conjecture} {\upshape \cite[Conjecture~E5.2]{AZ-Hilbert}}
This union of schemes is actually a projective scheme.
\end{conjecture}

The algebraic space of Theorem~\ref{th:Hilbert} is directly analogous
to the usual Hilbert scheme. However, \cite{AZ-Hilbert} also studies
``Hilbert schemes'' for abelian categories other than $\qgr R$.
For instance, they obtain the following theorem. Note that in this
case, the ``Hilbert scheme'' truly is a scheme.
\begin{theorem} {\upshape \cite[Theorem~E4.3]{AZ-Hilbert}}
Let $R$ be a connected graded, strongly noetherian ring. Let $P$ be
a finitely graded, noetherian $R$-module. 
The isomorphism classes of quotients of $P$ in $\gr R$ with fixed
Hilbert function $h$ are parameterized by a commutative projective scheme.
\end{theorem}

\section{Category theoretic schemes}

In Section~\ref{sec:AZT}, we defined a noncommutative projective scheme
as  $(\QGr R, \pi R)$ where $R$ was a ring with some nice properties.
Some work in noncommutative projective geometry is now moving into
a new paradigm, that of considering any category with ``geometric''
attributes as a noncommutative space. It is felt
that such a category should at least be a Grothendieck category
in the following sense.

\begin{definition}
An abelian category $\calC$ with generator and exact filtered direct limits
 is called a \emph{Grothendieck category}. 
 $\calC$ automatically has products and enough injectives.
\end{definition}

The following two definitions of a noncommutative surface are suggested in
\cite{SV}; we generalize the definition to that of a variety of dimension $n$.

\begin{definition}\label{def:nc1} 
{\upshape (Attempt I) } A \emph{noncommutative normal 
projective variety} (possibly singular)  of dimension $n$ is a category
$\calC$ of the form $\QGr R$, where $R$ is a noetherian connected
graded domain with $\GKdim R = n+1$ such that $R$ is a maximal order (this is
an analogue of normality).
$\calC$ is smooth if $\QGr R$ has homological dimension $n$
(that is, $n$ is minimal such that $\Ext^i(\calM, \calN) = 0$ for all $\calM, \calN \in \calC$
and all $i > n$).
\end{definition}

\begin{definition}\label{def:nc2}
 {\upshape (Attempt II) } A \emph{smooth noncommutative projective
variety} of dimension $n$
is a locally noetherian Grothendieck category $\calC$ of homological dimension
$n$ such that the full subcategory of noetherian objects is $\Ext$-finite and saturated. 
(Saturation differentiates between algebraic and analytic varieties; it turns out
that $\qgr R$ 
is saturated when $R$ and $R^\op$ are connected graded, noetherian,
 have $\chi$ and finite cohomological dimension,
and $\qgr R$ has finite homological dimension \cite{SV}).
\end{definition}

It is not clear which, if either, definition is ``correct.'' 
Work has been done on defining integral noncommutative schemes (that is,
integral Grothendieck categories) \cite{Smith-Integral}, so this may help
refine Definition~\ref{def:nc2}.

For both of these definitions, only the dimension $1$ case, the curve case,
has been settled. The noncommutative curves in the sense of
the first definition have been completely classified by \cite{AS} (even
when $R$ is not generated in degree one). Recently, curves in the sense
of the second definition have also been classified, falling into
only two different types. If one relaxes ``saturation'' to ``having
Serre duality via a Serre functor,'' then there are five different types
\cite{RV-hereditary}.

The surface case is a wide open problem for both definitions.
It is interesting to note that the two surface definitions conflict
slightly. More precisely, \cite{AV} gives an example of a twisted homogeneous
coordinate ring $R$ which is obtained from a smooth elliptic surface $X$,
yet $\GKdim R = 5$, not the expected $3$, so $\QGr R$ is not
a surface in the sense of Definition~\ref{def:nc1}. However, since $\QGr R \cong \Qch X$,
Definition~\ref{def:nc2} is satisfied.

There is now a theory of blowing up a point on a noncommutative smooth surface
in the sense of Definition~\ref{def:nc2} \cite{VdB2}. However, there
is not yet a ``Zariski's Main Theorem'' or an accepted concept of ``birational
equivalence class,'' 
so a classification of surfaces seems a long way off.
(See \cite[\S10.1]{SV} for a discussion of possible definitions for birational equivalence.)

Therefore, much work has been done on specific examples for both definitions.
We have already seen  a definition in Section~ \ref{sec:qpr} for noncommutative $\PP^2$
which is in line with Definition~\ref{def:nc1}.
Bondal and Polishchuk define more Grothendieck categories to be
 noncommutative $\PP^2$'s 
using a construction called ``$\ZZ$-algebras'' (which are actually
not rings in general).
\begin{comment}
\[
A_{ij} = \Hom(\OO(-j), \OO(-i))
\]
with ample sequence $(\OO(n))_n$.
\end{comment}

\begin{theorem} {\upshape \cite{BP-NCP2}}
Any Bondal-Polishchuk noncommutative $\PP^2$ ``comes from'' an AS-regular algebra $R$
of dimension $3$ with Hilbert series $(1-t)^{-3}$. That is, our original
Definition~\ref{def:qpr} of noncommutative $\PP^2$ and
that of Bondal and Polishchuk
 are in some sense the same.
\end{theorem}

A definition of noncommutative $\PP^1 \times \PP^1$ should include AS-regular algebras $R$
of dimension $3$ with Hilbert series $(1-t)^{-2}(1-t^2)^{-1}$ (that is, weight $(1,1,2)$).
But if one makes a Bondal-Polishchuk type definition of $\PP^1 \times \PP^1$ 
via $\ZZ$-algebras, one does get new categories
which don't come from such $R$. For more details on this,
as well as other recent surface generalizations, such as
$\PP^1$-bundles, cubic surfaces, and Cremona transformations, we
direct the reader to \cite{SV} again.

%%% \section* makes a section with no number
\section*{Acknowledgements}
The author would like to thank S.~P.~Smith and J.~T.~Stafford for their helpful suggestions.

\bibliographystyle{amsplain}
%\bibliography{ample}

%\begin{comment}
\providecommand{\bysame}{\leavevmode\hbox to3em{\hrulefill}\thinspace}
\providecommand{\MR}{\relax\ifhmode\unskip\space\fi MR }

%\end{comment}

\end{document}